\documentclass[english,a4paper,11pt]{article}

\usepackage[psamsfonts]{amssymb}
\usepackage{amsfonts,euscript,srcltx,url}
\usepackage{amsmath,amsthm,graphicx}
\usepackage{color,pgf,tikz}
\usetikzlibrary{calc,intersections, through,backgrounds}

\title{ A Self-Similar Dendrite with One-Point Intersection and\\ Infinite Post-Critical Set.}

\author{Prabhjot Singh \and Andrey Tetenov\footnote{Supported by Russian Foundation of Basic Research project 16-01-00414} }

\begin{document}
\newcommand{\qq}{\mathbb{Q}}
\newcommand{\rr}{\mathbb{R}}
\newcommand \nn {\mathbb{N}}
\newcommand \zz {\mathbb{Z}}
\newcommand \bbc {\mathbb{C}}
\newcommand \rd {\mathbb{R}^d}

 \newcommand {\al} {\alpha}
\newcommand {\be} {\beta}
\newcommand {\da} {\delta}
\newcommand {\Da} {\Delta}
\newcommand {\ga} {\gamma}
\newcommand {\ka} {\varkappa}
\newcommand {\Ga} {\Gamma}
\newcommand {\la} {\lambda}
\newcommand {\La} {\Lambda}
\newcommand{\om}{\omega}
\newcommand{\Om}{\Omega}
\newcommand {\sa} {\sigma}
\newcommand {\Sa} {\Sigma}
\newcommand {\te} {\theta}
\newcommand {\vte} {\vartheta}
\newcommand {\fy} {\varphi}
\newcommand {\Fy} { \Phi}
\newcommand {\ep} {\varepsilon}
\newcommand{\e}{\varepsilon}
\newcommand{ \vro}{\varrho}

\newcommand{\VEC}{\overrightarrow}
\newcommand{\IN}{{\subset}}
\newcommand{\NI}{{\supset}}
\newcommand \dd  {\partial}
\newcommand {\mmm}{{\setminus}}
\newcommand{\probel}{\vspace{.5cm}}
\newcommand{\8}{{\infty}}
\newcommand{\0}{{\varnothing}}
\newcommand{\vse}{$\blacksquare$}

\newcommand {\bfep} {{{\bar \varepsilon}}}
\newcommand {\Dl} {\Delta}
\newcommand{\vA}{{\vec {A}}}
\newcommand{\vB}{{\vec {B}}}
\newcommand{\vF}{{\vec {F}}}
\newcommand{\vf}{{\vec {f}}}
\newcommand{\vh}{{\vec {h}}}
\newcommand{\vJ}{{\vec {J}}}
\newcommand{\vK}{{\vec {K}}}
\newcommand{\vP}{{\vec {P}}}
\newcommand{\vX}{{\vec {X}}}
\newcommand{\vY}{{\vec {Y}}}
\newcommand{\vZ}{{\vec {Z}}}
\newcommand{\vx}{{\vec {x}}}
\newcommand{\va}{{\vec {a}}}
\newcommand{\vga}{{\vec {\gamma}}}

\newcommand{\eS}{{\EuScript S}}
\newcommand{\eH}{{\EuScript H}}
\newcommand{\eC}{{\EuScript C}}
\newcommand{\eP}{{\EuScript P}}
\newcommand{\eT}{{\EuScript T}}
\newcommand{\eG}{{\EuScript G}}
\newcommand{\eK}{{\EuScript K}}
\newcommand{\eF}{{\EuScript F}}
\newcommand{\eZ}{{\EuScript Z}}
\newcommand{\eL}{{\EuScript L}}
\newcommand{\eD}{{\EuScript D}}
\newcommand{\E}{{\EuScript E}}
\def \diam {\mathop{\rm diam}\nolimits}
\def \fix {\mathop{\rm fix}\nolimits}
\def \Lip {\mathop{\rm Lip}\nolimits}

\newcommand{\lf}{\lfloor}
\newcommand{\rc}{\rceil}

\newcommand{\io}{I^\infty}
\newcommand{\hio}{\hat{I}^\infty}
\newcommand{\imo}{I^{-\omega}}
\newcommand{\is}{I^*}

\newcommand{\Be}{{{\bf e}}}
\newcommand{\bi}{{{\bf i}}}
\newcommand{\bj}{{{\bf j}}}
\newcommand{\bk}{{{\bf k}}}
\newcommand{\bn}{{{\bf n}}}
\newcommand{\bx}{{{\bf x}}}
\newcommand{\by}{{{\bf y}}}

\newtheorem{thm}{\bf Theorem}
 \newtheorem{cor}[thm]{\bf Corollary}
 \newtheorem{lem}[thm]{\bf Lemma}
 \newtheorem{prop}[thm]{\bf Proposition}
 \newtheorem{dfn}[thm]{\bf Definition}

\newcommand{\dok}{{\bf{Proof}}}

\maketitle

\begin{abstract} We build an example of a system $\eS$ of similarities in $\rr^2$ whose attractor is a plane dendrite $K\supset [0,1]$ which satisfies one point intersection property, while the post-critical set of the system $\eS$ is a countable set  whose natural projection to $K$ is dense in the middle-third Cantor set.
\end{abstract}


\section{Introduction}
Let $\eS =\{S_1\ldots S_m\}$ be a system of contraction maps in $ \rr^n$.
 A non-empty compact set  $K$ $\in$ $ \rr^n$  satisfying
$K  =  S_1(K) \cup\cdots \cup S_m(K)$
is called invariant  set, or attractor, of the system $\eS $. The  uniqueness  and existence of  the  attractor $K $ is  provided by Hutchinson's  Theorem \cite{Hut}.\vspace{0.4em} \\  Let $I=\{1,2,\ldots,m\}$,  $I^*=\bigcup_{n=1}^\8I^n$ be the set of all finite $I$-tuples 
  and $I^{\8}=\{{\bf \al}=\al_1\al_2\ldots,\ \ \al_i\in I\}$ be the index  space and 
 $\pi:I^{\8}\rightarrow K$ be the address map.\vspace{0.4em}\\ 
A system $\eS$ satisfies open set condition(OSC) if there is nonempty open $O$ such that 
$S_i(O)\IN O$ and $S_i(O)\cap S_j(O)$ = $\0$ for any $i,j \in I, i\neq j$ \cite{Hut, Mor}. 
We say that the system $\eS$ satisfies one point  intersection  property \cite{BHR} if for any 
$i,j \in I, i\neq j$,   $\#(S_i(K)\cap S_j(K))\le 1$.\vspace{0.4em}\\
Let $\eC$ be the union of all  $S_i(K)\cap S_j(K)$, $i,j \in I, i\neq j$. 
The post-critical set $\eP$ of the system $\eS$ is the set of all
 $\alpha\in I^{\8}$ 
such that for some ${\bf j}\in I^*$,
 $S_ {\bf j}(\alpha)\in\eC$ . 
In other words, $\eP= \lbrace \sigma^k(\alpha) \vert   \alpha\in \eC,k\in \mathbb{N}\rbrace$, where the map $\sigma^k:I^{\8}\to I^{\8}$ is  defined by $\sigma^k(\al_1\al_2\ldots)=\al_{k+1}\al_{k+2}\ldots$.\vspace{0.4em}\\
A system $\eS$ is called {\em post-critically finite}(pcf) \cite{Kig} if its post-critical set is finite. This
 obviously implies  finite intersection property.
 
 Our aim is to show that the converse need not be true even in the  case of plane dendrites. We construct an example of pcf system $\eS$, whose attractor is a dendrite $K \IN\rr^2$, satisfying one point intersection property.\\

 \includegraphics[scale=.1]{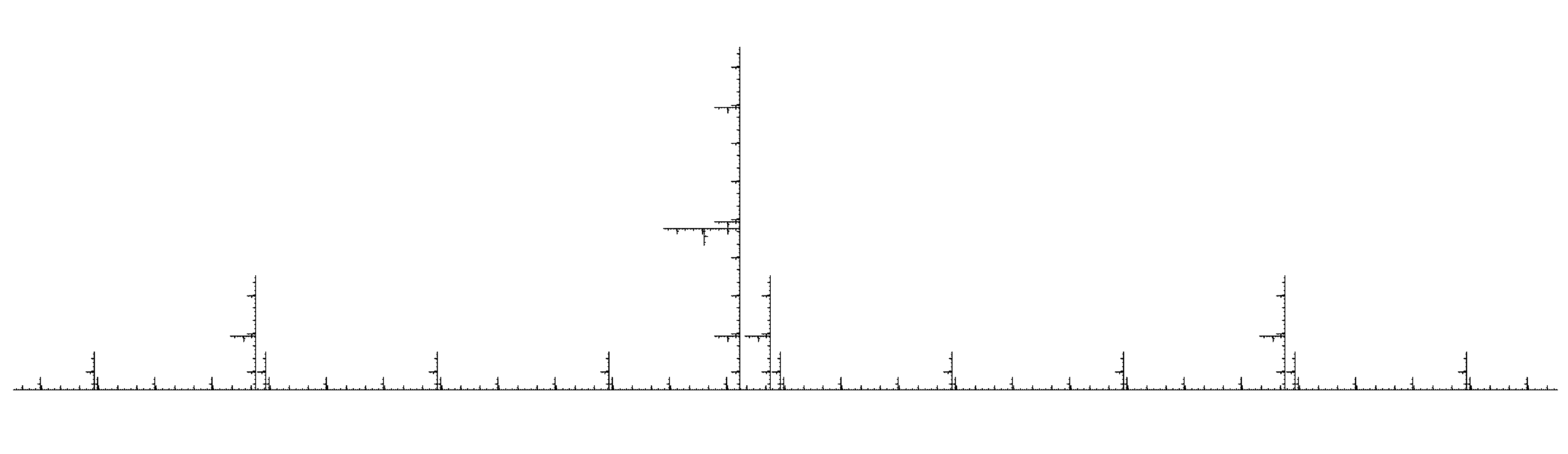}\\
 
 So we prove the following
 \begin{thm}
 There is a  system $\eS=\{S_1,S_2,S_3,S_h\}$ in $\rr^2$, whose attractor $K$ is a dendrite, which satisfies OSC and 1-point intersection property and has infinite  post-critical set    whose projection to $K$ is dense in the middle-third Cantor set. 
 
 \end{thm}
   
 \section{
 Construction}
 
Take a system $\eS=\lbrace S_0,S_1,S_2,S_h\rbrace$ of contraction similarities of $\rr^2$, defined by 
\begin{equation}
S_j(x,y)=((x+j)/3,y/3),\hspace{.2 cm} j=0,1,2 \mbox{    and   } 
S_h(x,y)= (-hy+c, hx)  
\end{equation}
and let $K$ be the attractor of $\eS$.  Here $c$ is   infinite base-3 fraction beginning with 0.11 and containing  all  finite tuples,  consisting of 0  and 2:\quad$c=0.110200220020222\ldots $

We will show that if $h$ is sufficiently small, then all the images\ \ $S_{i_1i_2\ldots i_nh}(K)$, $i_k=0,1,2$, are disjoint.\vspace{0.4em}\\
Put $I=\lbrace0,1,2\rbrace$ and denote by $I^*=\bigcup_{n=1}^\8 I^n$ the set of all tuples formed by $\lbrace0,1,2\rbrace$. Consider the images of $c$ under the maps $S_{\bj}, \bj\in I^*$. Using base-3 fractions, we can write  them as \quad
$S_{\bj}(c) = 3^{-n}c + {0.j_1\ldots j_n}, \ \   \bj\in I^* $, so ${(c_\bj)}_{k}=j_k$ for $k
\le n$ and ${(c_\bj)}_{k}=c_{k-n}$ for $k>n$.
\begin{center}
\begin{tikzpicture}[scale=11]
\draw [line width=.8pt] (0,0.)-- (1,0.);
\draw [line width=.75pt] (0.522,0.)-- (0.523,0.276)-- (0.45,0.2)-- (0.522,0.);
\draw [line width=.75pt] (0.467-0.021,0.)-- (0.469-0.021,0.082)-- (0.448-0.021,0.055)-- (0.467-0.021,0.);
\draw [line width=.75pt] (.1+.467,0.)-- (.1+0.469,0.082-.02)-- (.1+0.448,0.055-.02)-- (.1+0.467,0.);
\draw [line width=.75pt] (0.523,0.201)-- (0.45,0.2);
\draw[dotted] [line width=.75pt] (0.522,0.)-- (0.469-0.021,0.081);
\draw  (0.46,0.26) node {$ \Da$};
\draw  (0.425,0.2) node {$ B$};
\draw  (0.458,0.1) node {$ A_{\bj}$};
\draw  (0.5,0.22) node {$ h^2$};
\draw  (0.55,0.12) node {$ ch$};
\draw  (0.52,-0.025) node {$ c$};
\draw  (0.46,-0.025) node {$ c_{\bj}$};
\draw  (0.57,-0.025) node {$ c'_{\bj}$};
\draw  (0.40,0.07) node {$ \Da_{\bj}$};
\draw  (0.0,-0.025) node {$ (0,0)$};
\draw  (1.0,-0.025) node {$ (1,0)$};
\draw  (0.53,0.3) node {$\ (c,h)$};
\draw  (0.60,0.07) node {$ \Da'_{\bj}$};
\draw  (0.50,-0.09) node {Figure 1.};
\end{tikzpicture}
\end{center}
Let $ D$ be a triangle with vertices $\lbrace(0,0),(1,0),(c,h)\rbrace$ and $\Da=S_h(D)$.
 Since $c$ is not rational,  the point $c_{\bj}=S_{\bj}(c)$ is not equal to $c$ for any $\bj\in I^*$. 
  So either $c_{\bj}<c$ or $c_{\bj}>c$.\vspace{0.4em}\\
For $c_{\bj}<c$, $\Da_{\bj}\cap\Da\neq\0$ if  $A_{\bj}(c_j,h3^{-n})$ lies in $\Da$. To avoid this, the slope of the line $Bc$ has to be steeper than that of $A_{\bj}c$ (See Fig 1.):
$$
\frac{3^{-n}h}{(c-c_{\bj})} < \frac{ch}{h^2} \mbox{\quad or\quad} h^2 < 3^n(c-c_j)c \eqno{(1)}
$$
similarly, for $c<c'_{\bj}$, we have to require  that $S_{\bj}(B)$ does not lie in $\Da $
$$
\frac{h^2}{3^n} < (c'_{\bj}-c) \mbox{\quad or\quad} h^2 < 3^{n}(c'_{\bj}-c) \eqno{(2)}
$$

So we need to estimate $3^n(c-c_{\bj})$ and $3^n(c'_{\bj}-c)$. 
\vspace{0.4em}\\
\textbf{Case 1.} If $c_{\bj}<c$, there are the following possibilities: \medskip\\
(a) $j_1\ldots j_n$ = $i_1\ldots i_n$. Then  $(n+1)$-th entry $(c_{\bj})_{n+1}=i_1=1$ .
 Since  $i_{n+1}>1$, then $i_{n+1}=2$. 
 So $c_{\bj}<0.i_1\ldots i_n12$, $c>0.i_1\ldots i_n20$, then $c-c_{\bj} >{3^{-n-2}}$.\medskip \\
(b) $j_1\ldots j_k$ = $i_1\ldots i_k$ for some $k<n$ and $j_{k+1}<i_{k+1}$. Since the only entries allowed here are
  0 and 2, so $j_{k+1}=0$ and $i_{k+1}=2$.
  So $c> 0.i_1\ldots i_k2$ and $c_{\bj} < 0.i_1\ldots i_k1$,
   therefore $c-c_{\bj} < 3^{-k-1}.$\medskip\\

\textbf{Case 2.} $c<c'_{\bj}$, then \\
(a)  $j_1\ldots j_n=i_1\ldots i_n$. 
Since $(c_{\bj})_{n+1}=i_1=1$ ,  $c<c'_{\bj}$ 
implies $i_{n+1}=0$, so $c'_{\bj}>0.i_1\ldots i_n1102$, 
$c < 0.i_1\ldots i_n0(2)$, then $c-c'_{\bj} >3^{-n-2}$.\medskip\\
(b) $j_1\ldots j_k=i_1\ldots i_k$ for some $k<n$ and $j_{k+1}>i_{k+1}$. Then $j_{k+1}=2$, $i_{k+1}=0$, so $c'_{\bj}> 0.i_1\ldots i_k2$, $c < 0.i_1\ldots i_k1$, so $c'_j-c > 3^{-k-1}$.\vspace{0.4em}\\
(c) if $n=1$ and  $j_1=1$ i.e. $c_{\bj}=0.11102002\ldots$, then $c_{\bj}-c>0.0001121$. \\
Therefore $(c_{\bj}-c)3^n > 4/81.$\vspace{0.4em}\\
(d) if $n=2$ and $j_1j_2=11$, i.e. $c_{\bj}=0.11110200$, we similarly get $c_{\bj}-c>0.000201$. Thus, $(c_{\bj}-c)3^n > 2/9.$\vspace{0.4em}\\
So  if $h\le 2/9$ the inequalities (1) and (2) are satisfied. Further we  show that if $h\le 2/9$, the system $\eS$  satisfies open set condition(OSC) and one point intersection property and the attractor $K$ is a dendrite. 

\section{Proof of the Theorem.}

\begin{lem} If $h\le 2/9$, for any $\bi,\bj\in I^*$, $\Da_{\bi} \bigcap \Da_{\bj}=\0$.\end{lem}
\textbf{Proof.} Let $\bi=i_1 \ldots i_n$, $\bj=j_1\ldots j_m$. 
If $j_1\ldots j_k=i_1\ldots i_k$ for some $k<n$ and $j_{k+1}\neq i_{k+1}$, 
then $\Da_{\bi} \bigcap \Da_{\bj}=S_{i_1\ldots i_k}(\Da_{i_{k+1}\ldots i_{n}} \bigcap \Da_{j_{k+1}\ldots j_{m}})$. It follows from
 $\Da_{i_{k+1}\ldots i_n} \IN S_{i_{k+1}}( D)$ 
and  $\Da_{j_{k+1}\ldots j_m} \IN S_{j_{k+1}}( D)$ that $\Da_{i_{k+1}\ldots i_{n}}$ and $\Da_{j_{k+1}\ldots j_{m}}$ are disjoint.
 Thus $\Da_{\bi}\cap\Da_{\bj}=\0$.\\
If $m>n$ and $i_s=j_s$ for $s=1,\ldots,m$, then $\Da_{\bi}\bigcap \Da_{\bj}=S_{\bi}(\Da \bigcap \Da_{j_{n+1}\ldots j_{m}})$. 
By the construction, if $h\le  2/9$, $\Da \bigcap \Da_{j_{n+1}\ldots j_{m}}=\0$. $\blacksquare$
\begin{center}
\begin{tikzpicture}[x=9.0cm,y=9.0cm]
\draw (0.,0.)-- (0.9,0.);
\draw (0.17,0.05)-- (0.3,0.);
\draw (0.3,0.)-- (0.48,0.05);
\draw (0.48,0.05)-- (0.6,0.);
\draw (0.6,0.)-- (0.77,0.05);
\draw (0.9,0.)-- (0.51596,0.30049/2)-- (0.,0.);
\draw (0.51596,0.30049/2)-- (0.51258,0.)-- (0.4783,0.19926/2)-- (0.51596,0.30049/2);
\draw  (0.3,-0.04) node {$1/3$};
\draw  (0.6,-0.04) node {$2/3$};
\draw  (.6,0.15) node {$D$};
\draw  (.45,0.1) node {$\Da$};
\draw  (.17,0.024) node {$D_0$};
\draw  (.48,0.024) node {$D_1$};
\draw  (.77,0.024) node {$D_2$};
\draw  (0.45,-0.09) node {Fig. 2.};
\end{tikzpicture}
\end{center} 
\begin{lem} \textit{The system $\eS$ satisfies one point intersection property and open set condition(OSC)}.\end{lem}
\textbf{Proof.} Let $\Dot{D},\Dot{\Da}$ be the interiors of $D$ and $\Da$. 
Define $O=\Dot{\Da}\cup \bigcup_{\bi\in I^*} S_{\bi}(\Dot\Da)$. 
Obviously, for $i\in I$, $O_i=S_i(O)\IN O$. Moreover, $O_h=S_h(O)\IN\Dot\Da\IN O$.\\ 
Observe that with the only exception $S_1(\Dot D)\cap S_h(\Dot D)\neq\0$,
   the sets  $S_0(\Dot D), S_1(\Dot D)$, $S_2(\Dot D)$ and $S_h(\Dot D)$ are disjoint. 
Since $O\IN \Dot{D}$, the same is true for the sets $O_0,O_1,O_2$ and $O_h$. But $O_h\IN \Dot{\Da}$, so $O_h\bigcap O_1=\0$ too, therefore $O_0,O_1,O_2, O_h$ are disjoint and (OSC) is fulfilled.\\ 
 It follows from Lemma 2 that $\Da \bigcap S_1(\overline{O})=\lbrace c \rbrace$ and therefore $S_1(K) \bigcap S_h(K)=\lbrace c\rbrace$, which implies one point intersection property.$\blacksquare$
\begin{lem}   The system $\eS$ is post critically infinite and its post critical set is dense in the middle-third Cantor set $\mathcal{C}$.\end{lem}
\textbf{Proof.}  Let $y=.y_1y_2..., y_i\in\{0,2\}$ be base 3 representation for some point from the middle-third Cantor set $\mathcal{C}$.  Since the representation of $c$ contains all possible tuples of symbols 0  and 2, then for any $n\in \mathbb N$ there is $k=k(n)$ such that $c_{k+i}=y_{i}$ for $i=1,\ldots,n$. Therefore $\vert\sigma_{k}(c)-y\vert<3^{-n}$.  
So, the sequence $\sigma_{k(n)}(c)$, converges to the point $y\in \mathcal{C}$. $\blacksquare$

\medskip

To finish the proof of the  Theorem 1,  we need only to check  that the set $K$ is a dendrite.
 Let $\Da_{0}=\bigcup\limits_{\bi\in I^*} S_\bj(\Da)\cup\Da\cup[0,1]$.
  This set is compact and it is simply-connected, because the sets  $S_\bj(\Da)$ are disjoint.
   It is a strong  deformation retract of the  set $D$. 
   Define  $\Da_{k+1}=\bigcup\limits_{\bi\in I^*} S_\bj*S_h(\Da_k)\cup S_h(\Da_k)\cup[0,1]$.
The sets $\Da_k$ form a nested  sequence of compact simply-connected sets, each being a strong  deformation retract of the previous one. Then the intersection $\bigcap\limits_{k=1}^{\8} \Da_k=K$ is  a  strong  deformation retract of the set $D$. By Kigami's theorem \cite{Kig} it is locally connected  and  arcwise  connected. Since the interior of $K$ is empty, it contains no simple closed  curve, therefore it is a dendrite \cite[Theorem 1.1]{Char}.
\vse

\noindent
Prabhjot Singh \\ Central University of Rajasthan, 305801, Bandar-Sindri, India\\
\url{prabhjot198449@gmail.com}\\
Andrei Tetenov\\ Gorno-Altaisk state University, 649000,Gorno-Altaisk, Russia \\
\url{atet@mail.ru}\\

\end{document}